\input amstex
\documentstyle{amsppt}

\centerline {Introduction}

Let us consider the initial-boundary problem for the
Davey-Stewartson-II
system of equations:

$$
iu_{t} +
u_{xx} -u_{yy}
=\gamma|u|^{2}u +
\lambda u\varphi _{x} , \ \
\varphi _{xx}+ \varphi _{yy}
=\mu (\vert u \vert^{2})_{x},  \  (1)
$$

$$
u(x,y,t) \vert _{t=0} = u_{0}(x,y); \ \
\nabla \varphi \to 0,
\ \  x^2+y^2 \to \infty. \  (2)
$$

The Davey-Stewartson-II system of equations is the two-dimensional
generalization of the usual nonlinear Schrodinger equation.
Systems of the form (1) were first derived by Davey and Stewartson
in the shallow-water limit for the surface waves in fluid, [1].
Physical applications also
include plasma physics and nonlinear optics. In the above $ u(x,y,t)$
is the amplitude of a surface wave packet and $\varphi (x,y,t)$ is the
velocity potential of mean flow interacting with the surface wave.
It is worth mentioning that some of Davey-Stewartson systems of
equations are integrable, but many of them are not integrable.
The existence and uniqueness results with respect to the initial
data (in $ L_2, H^1, H^2 $) were recently obtained in the sufficiently
wide classes of generalized solutions, [2].
Our notice is centred on existence of classical solutions of the
initial-boundary problem for the more generally Davey-Stewartson-II
system of equations.

\centerline {The existence and uniqueness theorems}
Fourier transform with respect to spatial variables lies
in the base of proof of solvability of the initial-boundary
problem (1)-(2). By using this traditional approach to the
proof of existence theorems, we may pass from the system of partial
differential equations (1) to a single integral equation in terms
of the Fourier transform $\hat u(k,m,t)$ of $u(x,y,t)$:

$$
\hat u(k,m,t)=\hat u_{0}(k,m)e^{i(m^2-k^2)t}
-i\int_{0}^{t}e^{i(k^2-m^2)(\tau -t)}
L(\hat u(k,m, \tau ) e^{2i(k^2-m^2) \tau} )d \tau. \ (3)
$$

where the operator $L$ is defined by the following formula:

$$
L(\hat u)  =\gamma
\hat u\star \hat u^{*}\star \hat u+
\lambda\mu\hat u \star [\frac {k^2}{k^2+m^2}
(\hat u\star \hat u^{*})].
$$

Here the sign $\star$ is a convolution with respect to $(k,m)$.

Theorem 1. Let  $u_0(x,y) \in H^p$ provided $ p>1. $
There exists a positive $T>0$ such that
$ \forall t \in [0,T]$ the initial-boundary problem (1)-(2)
has a unique solution such that: \ \
$
u, \nabla \varphi \in C^1([0,T];H^p).
$

The basic idea of the proof is now to use a contraction mapping
method in a suitable scale of Banach spaces. In these Banach
spaces integral operator from (3) should be bounded and Lipschitzian.
By analogy with [3,4], we introduce the scale of Banach spaces
$C( [0,T];H^p )$
of Fourier transforms $\hat u(k,m,t)$, which are continuous with respect
to $(k,m),$ $ t\in [0,T]$ and decreasing at the infinity,
with the following norms:
$$
\parallel \hat v(k,m,t) \parallel = \sup\limits_{t\in [0,T]}\parallel
(1+\vert k \vert + \vert m \vert)^p \hat v(k,m,t) \parallel_{L_{2}} .
$$

In view of nonlinearity of $L$ operator, our Banach space must be
an algebra with respect to the convolution operation. This problem
was solved provided $p>1$. It turns out that in this case the convolution
is a Lipschitzian operator in Sobolev space $H^p$. Then
if we choose a suitable $T$ such that the integral operator from (3)
would be contraction, the solvability of the equation (3), and so of the
problem (1)-(2) follows from the Ovsjannikov's
results for the equations in Banach spaces, [5].

This result is directly applicable to the more
general systems of equations in multidimensions. Namely,
we consider the following initial-boundary problem for the
system of partial differential equations:

$$
i\partial_{t} u +
\omega(-i\partial _{x_{1}},\dots, -i\partial _{x_{n}}) u
=g(u,u^*) +
\bar {h}(u,u^*) \nabla\varphi , \ \
P\varphi
=\nabla \bar {f} (u,u^*), \ (4)
$$

$$
u(x_{1},\dots,x_{n},t)|_{t=0} = u_{0}(x_{1},\dots,x_{n}); \ \
\nabla \varphi \to 0,
\quad x_{1}^{2}+\dots +x_{n}^{2} \to \infty, \ (5)
$$

\noindent
where $P$ is a second-order differential operator with the constant
coefficients and $\omega (k_{1},\dots,
k_{n})$ is a polinom with respect to $ k_{1},\dots,k_{n}.$

Theorem 2. Let
1) functions $\bar {f},g,\bar {h} $ are analytical
with respect to all arguments such that:
$\bar {f}(0,0)=g(0,0)=\bar {h}(0,0)=0;$
2) initial function $u_0(x_{1},\dots,x_{n}) \in H^p$ provided $p>n/2$.
There exists a positive $T>0$ such that
$ \forall t \in [0,T]$ the initial-boundary problem (4)-(5)
has a unique solution such that:
$
u,\nabla \varphi (x_{1},\dots,x_{n},t) \in C^1([0,T];H^{p}).
$

\centerline { Literature }

[1] Davey A., Stewartson K.,
On the three-dimensional packets of surfase waves,
Proc.R.Soc.London, 1974, Vol.338, Ser.A, pp.101--110.

[2] Ghidaglia J.-M., Saut J.-C.,
On the initial value problem for the Davey-Stewartson systems
Nonlinearity, 1990, Vol.3, pp.475--506.

[3] Vakulenko S.A.,
Justification of asymptotic formula for the perturbed
Klein-Fock-Gordon equation,
Zap.nauch.semin. LOMI, 1981, Vol.104, pp.84--92.

[4] Kalyakin L.A.,
Asymptotic collapse of one-dimensional wave packet in
nonlinear dispersive medium,
Mat.sbornik, 1987, Vol.132, No.4, pp.470--495.

[5] Ovsjannikov L.V.
A nonlinear Cauchy problem in a scale of Banach spaces,
Dokl.Akad.Nauk SSSR, 1971, Vol.200, No.4, pp.789--792.

\end